\documentclass{amsart}

\usepackage{amssymb}
\theoremstyle{plain}

\newtheorem{thm}{Theorem}
\newtheorem{prop}[thm]{Proposition}

\newtheorem{lem}[thm]{Lemma}

\theoremstyle{definition}

\newtheorem{calc}[thm]{Calculation}

\author{Michael Robinson}
\address{Center for Applied Mathematics\\657 Rhodes Hall\\Cornell
  University, Ithaca, NY 14850}
\email{robinm@cam.cornell.edu}



\title[Eternal solutions for parabolic equations]{Construction of eternal solutions for a semilinear parabolic equation}

\begin{document}

\begin{abstract}
Eternal solutions of parabolic equations (those which are defined for
all time) are typically rather rare.  For example, the heat equation
has exactly one eternal solution -- the trivial solution.  While
solutions to the heat equation exist for all forward time, they cannot
be extended backwards in time.  Nonlinearities exasperate the
situation somewhat, in that solutions may form singularities in both
backward and forward time.  However, semilinear parabolic equations
can also support nontrivial eternal solutions.  This article shows how
nontrivial eternal solutions can be constructed for a semilinear
equation that has at least two distinct equilibrium solutions.  The
resulting eternal solution is a heteroclinic orbit which connects the
two given equilibria.
\end{abstract}

\maketitle

\section{Introduction}

Consider
\begin{equation}
\label{pde}
\frac{\partial u(t,x)}{\partial t}=\Delta u(t,x) +  \sum_{i=0}^N
a_i(x)u^i(t,x),
\end{equation}
where $t \in \mathbb{R}$ and $x \in \mathbb{R}^n$, and the $a_i$ are
bounded and smooth.  In this article, we consider {\it eternal}
solutions, those classical solutions $u$ which satisfy \eqref{pde} for all time $t \in \mathbb{R}$. 

This kind of equation provides a simple model for a number of physical
phenomena.  First, choosing the right side to be $\Delta u - u^2 + a_1
u$ results in an equation which can represent a model of the
population of a single species with diffusion and a spatially-varying
carrying capacity, $a_1(x)$.  As a second application, this equation
is a very simple model of combustion.  If $a_1$ is a positive
constant, then the equation supports traveling waves.  Such traveling
waves can model the propagation of a flame through a fuel source.

Equations of the form \eqref{pde} have been of interest to researchers
for quite some time.  Existence and uniqueness of solutions on short
time intervals (on strips $(-t_0,t_0)\times\mathbb{R}^n$) can been
shown using semigroup methods and are entirely standard
\cite{ZeidlerIIA}. However, there are obstructions to the existence of
eternal solutions.  Aside from the typical loss of regularity due to
solving the backwards heat equation, there is also a blow-up
phenomenon which can spoil existence in the forward-time solution to
\eqref{pde}.  Blow-up phenonmena in the forward time Cauchy problem
(where one does not consider $t<0$) have been studied by a number of
authors \cite{Fujita} \cite{Kobayashi_1977} \cite{Weissler_1981}
\cite{Klainerman_1982} \cite{Brezis_1984} \cite{Zheng_1986}
\cite{Zheng_1995}.  More recently, Zhang {\it et al.}
(\cite{Zhang_2000} \cite{Souplet_2002} \cite{Wrkich_2007}) studied
global existence for the forward Cauchy problem for
\begin{equation*}
\frac{\partial u}{\partial t} = \Delta u + u^p - V(x) u
\end{equation*}
for positive $u,V$.  Du and Ma studied a related problem in
\cite{DuMa2001} under more restricted conditions on the coefficients
but they obtained stronger existence results.  In fact, they found
that all of the solutions which were defined for all $t>0$ tended to
equilibrium solutions.

Eternal solutions to \eqref{pde} are rather rare.  Most works which
describe blow-up make the assumption that the solution is positive.
Unfortunately, blow-up is much more difficult to characterize in the
general situation, and understanding exactly what kind of initial
conditions are responsible for blow-up in the Cauchy problem for
\eqref{pde} is an important part of the question.

As an aside, the boundary value problem that results from taking $x
\in \Omega \subset \mathbb{R}^n$ for some bounded $\Omega$ (instead of
$x \in \mathbb{R}^n$) has also been discussed extensively in the literature
\cite{Henry} \cite{Jost_2007} \cite{BrunovskyFiedler1989}.  For the
boundary value problem, all bounded forward Cauchy problem solutions
tend to limits as $|t|\to\infty$, and these limits are equilibrium
solutions.

The existence of eternal solutions is a difficult problem, because the
backward-time Cauchy problem is well known to be ill-posed.
Obviously, equilibrium solutions are trivial examples of such eternal
solutions, and in \cite{RobinsonNonauto} it was shown that they can
exist.  It is not at all clear that there are other eternal solutions,
and indeed there may not be.  In this article we assume the existence
of a pair of nonintersecting equilibrium solutions and construct a
heteroclinic orbit which connects them.  (A {\it heteroclinic orbit}
is a special kind of eternal solution, whose limits as $t\to\pm\infty$
are equilibria.)

For simplicity and concreteness, we will work with the more
limited equation
\begin{equation}
\label{pde1}
\frac{\partial u(t,x)}{\partial t} = \frac{\partial^2 u(t,x)}{\partial
  x^2} - u^2(t,x) + \phi(x),
\end{equation}
where $\phi$ is a smooth function which decays to zero.  We follow the
general technique for constructing ``ancient solutions,'' which was
used in a different context by Hamilton and Perelman.  It should be
emphasized that the technique examined in this article can handle the
problem for \eqref{pde} in full generality (under mild decay
assumptions for the coefficients $a_i$), though this complicates the
exposition needlessly.

This simpler model still provides insight into applications, as it is
still a model of the population of a single species, with a
spatially-varying carrying capacity, $\phi$.  Indeed, one easily finds
that under certain conditions the behavior of solutions to
\eqref{pde1} is reminiscent of the growth and (admittedly tenuous)
control of invasive species \cite{Blaustein_2001}.  It is the control
of invasive species that is of most interest, and it is also what the
structure of the space of heteroclines describes.  In one of the
examples given in \cite{RobinsonNonauto}, there is one ``more stable''
equilibrium, and several other ``less stable'' ones.  The more stable
equilibrium can be thought of as the situation where an invasive
species dominates.  The task, then, is to try to perturb the system so
that it no longer is attracted to that equilibrium.  An optimal
control approach is to perturb the system so that it barely crosses
the boundary of the stable manifold of the the undesired equilibrium,
and thereby the invasive species is eventually brought under control
with minimal disturbance to the rest of the environment.  Such an
optimal control approach, though, is beyond the scope of this article.

\section{Equilibrium solutions}

We choose $\phi(x)=(x^2-0.4)e^{-x^2/2}$.  It has been shown in \cite{RobinsonNonauto}, that in this situation, there exists a pair of equilibrium solutions $f_+,f_-$ with the following properties:

\begin{enumerate}
\item $f_+$ and $f_-$ are smooth and bounded,
\item $f_+$ and $f_-$ have bounded first and second derivatives,
\item $f_+$ and $f_-$ are asymptotic to $6/x^2$ for large $x$, and so
  both belong to $L^1(\mathbb{R})$,
\item $f_+(x)>f_-(x)$ for all $x$,
\item there is no equilibrium solution $f_2$ with
  $f_+(x)>f_2(x)>f_-(x)$ for all $x$,
\end{enumerate}
and additionally, there exists a one-parameter family $g_c$ of
solutions to 
\begin{equation}
\label{gc_eqn}
0=g_c''(x) - g_c^2(x) + \phi_c(x)
\end{equation}
with 
\begin{enumerate}
\item $c\in[0,1)$,
\item $g_0=f_-$ and $\phi_0=\phi$,
\item $\phi_a(x) < \phi_b(x)$ and $g_a(x) > g_b(x)$ for all $x$ if $a>b$.
\end{enumerate}

The latter set of properties can occur as a consequence of the
specific structure of $f_-$.  For instance, consider the following result.

\begin{prop}
  Suppose $f_- \in C^{2,\alpha}(\mathbb{R})$ satisfies the above
  conditions and additionally, there is a compact $K \subset
  \mathbb{R}$ with nonempty interior such that $f_-$ is negative on
  the interior of $K$ and is nonnegative on the complement of $K$.
  Then such a family $g_c$ above exists.
\begin{proof}(Sketch)
Work in $T_{f_-}C^{2,\alpha}(\mathbb{R})$, the tangent space at
$f_-$.  Then \eqref{gc_eqn} becomes its linearization (for $h_c$, say), namely
\begin{equation}
\label{gc_lin_eqn}
0=h_c''(x) - 2 f_-(x) h_c(x) + (\phi_c - \phi).
\end{equation}
Consider the slighly different problem, 
\begin{equation}
\label{gc_lin_eg_eqn}
0=y''(x)-2f_-(x) y(x) + v(x) y(x),
\end{equation}
where $v$ is a smooth function to be determined.  If we can find a $v
\le 0$ such that $y>0$ and $y \to 0$ as $|x| \to \infty$, then we are
done, because we simply let $vy = \phi_c - \phi$ in
\eqref{gc_lin_eqn}.  In that case, $h_c=y$ has the required
properties.  
We sketch why such a $v$ exists:
\begin{itemize}
\item If $v \equiv 0$, then $y \equiv 0$ is a solution, giving
  $g_c=f_-$ as a base case.
\item If $v(x) = - 2 \|u\|_\infty \beta(x)$ for $\beta$ is a smooth
  bump function with compact support and $\beta | K = 1$, then the
  Sturm-Liouville comparison theorem implies that $y$ has no sign
  changes.  We can take $y$ strictly positive.  However, in this case,
  the Sturm-Liouville theorem imples that there are no critical points
  of $y$ either, so $y$ may not tend to zero as $|x|\to \infty$.
\item Hence there should exist an $s$ with $0<s<2\|u\|_\infty$ such
  that if $v(x) = -s \beta(x)$, then $y$ has no sign changes, one
  critical point, and tends to zero as $|x|\to\infty$.  This choice of
  $v$ is what is required.  (The precise details of this argument fall
  under standard Sturm-Liouville theory, which are omitted here.)
\end{itemize}
\end{proof}
\end{prop}

In what follows, we shall not be concerned with the exact form of
$\phi$, but rather we shall assume that the above properties of the
equilibria hold.  Many other choices of $\phi$ will allow a similar
construction.

\begin{lem}
\label{lucid_lem}
The set
\begin{equation}
\label{funnel_eqn}
W=\{v \in C^2(\mathbb{R}) | f_-(x) < v(x) < f_+(x)\text{ for all }x\}
\end{equation}
is a forward invariant set for \eqref{pde1}.  That is, if $u$ is a solution to
\eqref{pde1} and $u(t_0) \in W$, then $u(t) \in W$ for all $t>t_0$.
\begin{proof}
We show that the flow of \eqref{pde1} is inward whenever a timeslice is tangent to either $f_-$ or $f_+$.  To this end, define the
set $B$
\begin{eqnarray*}
B&=&\{v \in C^2(\mathbb{R}) | f_-(x) \le v(x) \le f_+(x)\text{  for all }x,\text{ and there exists an }x_0\\&&\text{ such that }
  v(x_0)=f_+(x)\text{ or } v(x_0)=f_-(x)\}.
\end{eqnarray*}
Without loss of generality, consider a $v\in B$ with a single
point of tangency, $v(x_0)=f_-(x_0)$.  At such a point $x_0$, the
smoothness of $v$ and $f_-$ implies that $\Delta v(x_0) \ge \Delta
f_-(x_0)$ using the maximum principle.  Then, if $u$ is a solution to
\eqref{pde1} with $u(0,x)=v(x)$, we have that
\begin{eqnarray*}
\frac{\partial u(0,x_0)}{\partial t} &=& \Delta v(x_0) - v^2(x_0) +
\phi(x_0)\\ 
&\ge&\Delta f_-(x_0) - f_-^2(x_0) + \phi(x_0)=0,\\
\end{eqnarray*}
hence the flow is inward.  One can repeat the above argument for each
point of tangency, and for tangency with $f_+$ as well.
\end{proof}
\end{lem}

\begin{lem}
\label{nested_funnel_lem}
Solutions to the Cauchy problem 
\begin{equation}
\label{pde_cauchy}
\begin{cases}
\frac{\partial u(t,x)}{\partial t} = \frac{\partial^2 u(t,x)}{\partial
  x^2} - u^2(t,x) + \phi(x),\\
u(0,x) = U(x) \in W_c
\end{cases}
\end{equation}
where
\begin{equation*}
W_c=\{v \in C^2(\mathbb{R}) | g_c(x) < v(x) < f_+(x)\text{ for all }x\}
\end{equation*}
for $c \in [0,1)$ have the property that they lie in $L^1 \cap
L^\infty(\mathbb{R})$ for all $t>0$.  We shall assume that $U$ has
bounded first and second derivatives.

Additionally, when $c \in (0,1)$, solutions to \eqref{pde_cauchy}
cannot have $f_-$ as a limit as $t\to\infty$.
\begin{proof}
  The fact that solutions lie in $L^1 \cap L^\infty(\mathbb{R})$ is
  immediate from Lemma \ref{lucid_lem} and the asymptotic behavior of
  $f_+,f_-$ (Section 4 of \cite{RobinsonNonauto}).  Observe that for each $c\in
  [0,1)$, $W_c$ is forward invariant, and that $W_a \subset W_b$ if
  $a>b$.  Since $f_-$ is not in $W_c$ for $c$ strictly larger than 0,
  the proof is completed.
\end{proof}
\end{lem}

The following is an outline for the rest of the article.  We show that
all solutions to \eqref{pde_cauchy} have bounded first and second
spatial derivatives.  This implies that all of their first partial
derivatives are bounded (the time derivative is controlled by
\eqref{pde1}).  Using the fact that \eqref{pde1} is autonomous in
time, time translations of solutions are also solutions.  We therefore
construct a sequence of solutions $\{u_k\}$ to Cauchy problems started
at $t=0,T_{1}, T_{2}, ...$ which tend to $f_+$ as $t\to +\infty$, but
their initial conditions tend to $f_-$ as $k\to \infty$.  By Ascoli's
theorem, this sequence converges uniformly on compact subsets to a
continuous eternal solution.

\section{Integral equation formulation}
In order to estimate the derivatives of a solution to
\eqref{pde_cauchy}, it is more convenient to work with an integral
equation formulation of \eqref{pde_cauchy}.  This is obtained in the
usual way.

\begin{eqnarray*}
\frac{\partial u}{\partial t} &=& \Delta u - u^2 + \phi\\
\left(\frac{\partial}{\partial t} - \Delta \right) u &=& - u^2 + \phi\\
u &=& \left(\frac{\partial}{\partial t} - \Delta \right)^{-1} ( \phi - u^2)\\
\end{eqnarray*}
\begin{equation}
\label{int_eqn}
u(t,x)=\int_{-\infty}^\infty H(t,x-y) U(y) dy + \int_0^t
\int_{-\infty}^\infty H(t-s,x-y)\left(\phi(y) - u^2(s,y)\right) dy\, ds,
\end{equation}
where $H(t,x)=\frac{1}{\sqrt{4 \pi t}}e^{- \frac{x^2}{4t}}$ is the
usual heat kernel.  

\begin{calc}
\label{first_deriv_short}
We begin by estimating the first derivative of $u$ for a short time.
Let $T>0$ be given, and consider $0\le t\le T$.  The key fact is that
$\int H(t,x) dx = 1$ for all $t$.  Using \eqref{int_eqn}

\begin{eqnarray*}
\left \| \frac{\partial u}{ \partial x} \right \|_\infty &\le& \left \|
\frac{\partial U}{\partial x} \right \|_\infty + \left|\int_0^t \int_{-\infty}^\infty
\frac{\partial}{\partial x}H(t-s,x-y)\left(\phi(y) - u^2(s,y)\right)
dy\, ds\right|\\
&\le& \left \|\frac{\partial U}{\partial x} \right \|_\infty +
\int_0^t \int_{-\infty}^\infty
\left|\frac{\partial}{\partial y}(H(t-s,x-y))\left(\phi(y) -
u^2(s,y)\right) \right| dy\, ds\\
&\le& \left \|\frac{\partial U}{\partial x} \right \|_\infty +
\int_0^t \int_{-\infty}^\infty
\left|H(t-s,x-y)\left(\frac{\partial \phi}{\partial y} - 2u
\frac{\partial u}{\partial y}\right) \right| dy\, ds\\
&\le& \left \|\frac{\partial U}{\partial x} \right \|_\infty +
T \left\|\frac{\partial \phi}{\partial x}\right\|_\infty +
2 \| u\|_\infty \int_0^t \left\|
\frac{\partial u}{\partial x} \right\|_\infty ds.\\
\end{eqnarray*}
This integral equation fence is easily solved to give
\begin{eqnarray*}
\left \| \frac{\partial u}{ \partial x} \right \|_\infty &\le& \left (
\left\| \frac{\partial U}{\partial x}\right \|_\infty + T\left\|
\frac{\phi}{\partial x}\right \|_\infty \right) e^{2 t
  \max\{\|f_+\|_\infty,\|f_-\|_\infty\}}\\
&\le&K_1 e^{K_2 T}.
\end{eqnarray*}
\end{calc}

\begin{calc}
\label{second_deriv_short}
With the same choice of $T$ as above, we find a bound for the second
derivative in the same way:

\begin{eqnarray*}
\left \| \frac{\partial^2 u}{\partial x^2} \right \|_\infty &\le&
\left \|\frac{\partial^2 U}{\partial x^2} \right \|_\infty + T \left \|
\frac{\partial^2 \phi}{\partial x^2}\right \|_\infty + \int_0^t \left
\| \frac{\partial}{\partial y} \left( 2 u \frac{\partial u}{\partial
  y}\right)\right\|_\infty ds\\
 &\le&
\left \|\frac{\partial^2 U}{\partial x^2} \right \|_\infty + T \left \|
\frac{\partial^2 \phi}{\partial x^2}\right \|_\infty + \int_0^t 2
\left \| \frac{\partial u}{\partial x} \right \|^2_\infty + 2 \|
u \|_\infty \left \| \frac{\partial^2 u}{\partial x^2}\right \|_\infty
ds\\
&\le& K_3 e^{K_2 T}
\end{eqnarray*}
for some $K_3$ which depends on $U$, $\phi$, and $T$.
\end{calc}

\begin{calc}
\label{second_deriv_long}
Now, we extend Calculation \ref{second_deriv_short} to handle $t>T$,

\begin{eqnarray*}
\left \| \frac{\partial^2 u}{\partial x^2} \right \|_\infty &\le&
\left \| \frac{\partial^2 U}{\partial x^2} \right \|_\infty +
\left|\frac{\partial^2}{\partial x^2}\int_0^T \int_{-\infty}^\infty
H(t-s,x-y) (\phi(y) - u^2(s,y)) dy\, ds\right | +\\&& 
\left | \frac{\partial^2}{\partial x^2}\int_T^t \int_{-\infty}^\infty
H(t-s,x-y) (\phi(y) - u^2(s,y)) dy\, ds \right |
\\
&\le&
K_3 e^{K_2 T} + 
\int_T^t \left\|\frac{\partial^2}{\partial x^2} H(t-s,x)\right\|_\infty
(\|\phi\|_1 + \|u^2\|_1) ds \\
&\le&
K_3 e^{K_2 T} + K_4 \int_T^t \frac{1}{s\sqrt{s}} ds +K'_4 \int_T^t \frac{1}{s^2\sqrt{s}} ds\\
&\le&
K_3 e^{K_2 T} + K_5 \left(\frac{1}{\sqrt{T}} -
\frac{1}{\sqrt{t}}\right) +K'_5 \left(\frac{1}{T\sqrt{T}} - \frac{1}{t\sqrt{t}}\right)\\
&\le&
K_3 e^{K_2 T} + K_6,\\
\end{eqnarray*}
hence there is a uniform upper bound on $\left \| \frac{\partial^2
  u}{\partial x^2} \right \|_\infty$ which depends only on the initial
  conditions, $\phi$, and $T$.
\end{calc}

\begin{lem}
\label{bdd_deriv}
Let $f\in C^2(\mathbb{R})$ be a bounded function with a bounded second
derivative.  Then the first derivative of $f$ is also bounded, and the
bound depends only on $\|f\|_\infty$ and $\|f''\|_\infty$.
\begin{proof}
The proof is elementary.  The key fact is that
at its maxima and minima, $f$ has a horizontal tangent.  From a
horizontal tangent, the quickest $f'$ can grow is at a rate of
$\|f''\|_\infty$.  However, since $f$ is bounded, there is a maximum
amount that this growth of $f'$ can accrue.  Indeed, a sharp estimate
is
\begin{eqnarray*}
\|f'\|_\infty \le \sqrt{2 \|f\|_\infty \|f''\|_\infty}.
\end{eqnarray*}
\end{proof}
\end{lem}

Using the fact that $u$ is bounded, Lemma \ref{bdd_deriv} implies that
the first spatial derivative of $u$ is bounded.  By \eqref{pde1}, it
is clear that the first time derivative of $u$ is also bounded.

\begin{lem}
\label{bdd_action}
As an immediate consequence of Lemmas \ref{nested_funnel_lem} and
\ref{bdd_deriv}, the action integral
\begin{equation*}
A(u(t))=\int_{-\infty}^\infty \frac{1}{2}\left| \frac{\partial
  u}{\partial x} \right|^2 + \frac{1}{3}u^3(t,x) - u(t,x)\phi(x) dx
\end{equation*}
is bounded.  Therefore, the solutions to the Cauchy problem
\eqref{pde_cauchy} all tend to limits as $t\to\infty$ (Corollary 7 of \cite{RobinsonClassify}).  By Lemma \ref{nested_funnel_lem}, we
conclude that they all tend to the common limit of $f_+$ when $c>0$.
\begin{proof}
  The latter two terms are bounded due to the fact that $u$ lies in
  $L^1 \cap L^\infty(\mathbb{R})$ for all $t$.  The bound on the first
  term comes from combining the fact that $u$ and its first two
  spatial derivatives are bounded with the asymptotic decay of
  $f_\pm$, and is otherwise straightforward (use L'H\^opital's rule).
\end{proof}
\end{lem}

\section{Construction of an eternal solution}

Let 
\begin{equation*}
U_k(x) = (1-2^{-k-1})g_{1/k}(x) + 2^{-k-1}f_+(x),\text{ for }k \ge 0
\end{equation*}
noting that $U_k \to f_-$ as $k\to \infty$.  Since $U_k$ is a convex
combination of $f_+$ and $g_{1/k}$, it follows that $U_k \in W_{1/k}$ for all
$k$.  Also, since $f_+$ and $f_-$ have bounded first and second
derivatives, the $\{U_k\}$ have a common bound for their first and second
derivatives.

Now consider solutions to the following set of Cauchy problems
\begin{equation}
\label{cauchy_set}
\begin{cases}
\frac{\partial u_k(t,x)}{\partial t} = \frac{\partial^2 u_k(t,x)}{\partial
  x^2} - u_k^2(t,x) + \phi(x),\\
u_k(T_k,x) = U_k(x).
\end{cases}
\end{equation}  

We choose $T_k$ so that for all $k>0$, $u_k(0,0)=u_0(0,0)$.  We can do
this using the continuity of the solution and Lemma \ref{bdd_action}.  As
$k\to \infty$, solutions are started nearer and nearer to the
equilibrium $f_-$, so we are forced to choose $T_k \to -\infty$
as $k\to \infty$.

It's clear that each solution $u_k$ is defined for only $t>T_k$.  
However, for each compact set $S \subset \mathbb{R}^2$, there are
infinitely many elements of $\{u_k\}$ which are defined on it.  The
results of the previous section imply that $\{u_k\}$ is a bounded,
equicontinous family.  As a result, Ascoli's theorem implies that
$\{u_k\}$ converges uniformly on compact subsets to a continous $u$,
which is an eternal solution to \eqref{pde1}.

Our constructed eternal solution will have the value $u(0,0)=u_0(0,0)$,
which is strictly between $f_+$ and $f_-$.  As a result, the eternal
solution we have constructed is not an equilibrium solution.  By Lemma
\ref{bdd_action}, it is a finite energy solution, so it must be a
heteroclinic orbit connecting $f_-$ to $f_+$.

\bibliography{global_bib}
\bibliographystyle{plain}

\end{document}